\documentclass[12pt]{amsart}

\usepackage[T1]{fontenc}
\usepackage[utf8]{inputenc}
\usepackage{lmodern}
\usepackage{microtype}
\usepackage{mathtools}
\usepackage{amssymb}
\usepackage{parskip}
\usepackage{geometry}
\newtheorem{theorem}{Theorem}

\newtheorem{lemma}[theorem]{Lemma}
\theoremstyle{definition}

\theoremstyle{remark}

\newcommand{\Ztwo}{\mathbb Z/2\mathbb Z}
\newcommand{\Nerve}{\mathcal N}
\newcommand{\MaxNerve}{\mathcal N^{\max}}
\newcommand{\Cosets}{\mathcal C}
\newcommand{\ab}{\mathrm{ab}}

\title[Rigidity of Right-Angled Coxeter Groups]{Rigidity of Right-Angled Coxeter Groups}
\author{David G. Radcliffe}
\email{dradcliffe@gmail.com}
\date{11 January 1999, Revised 12 May 2026}

\begin{document}

\begin{abstract}
	Let $W$ be a right-angled Coxeter group. If $S$ and $S'$ are finite Coxeter generating sets for $W$, then the Coxeter systems $(W,S)$ and $(W,S')$ are equivalent.
\end{abstract}

\maketitle

\section{Introduction}

A \emph{Coxeter system} is a pair $(W,S)$ in which $W$ has a presentation
\[
	W=\left\langle S \;\middle|\; (st)^{m(s,t)}=1 \text{ for } s,t\in S \right\rangle,
\]
where $S$ is finite, $m(s,s)=1$ for every $s\in S$, and, for $s\ne t$,
\[
	m(s,t)=m(t,s)\in \{2,3,4,\ldots,\infty\}.
\]
The relation $(st)^{m(s,t)}=1$ is omitted when $m(s,t)=\infty$. In a Coxeter system the number $m(s,t)$ is the order of $st$ in $W$; see \cite[p.~110]{Humphreys}.

A set $S$ occurring in such a presentation is called a \emph{Coxeter generating set} for $W$. Two Coxeter systems $(W,S)$ and $(W',S')$ are \emph{equivalent} if there is an isomorphism $\phi\colon W\to W'$ such that $\phi(S)=S'$.

The Coxeter system $(W,S)$ is \emph{right-angled} if $m(s,t)\in\{2,\infty\}$ whenever $s\ne t$. Equivalently, the only defining relations are $s^2=1$ for $s\in S$ and commutation relations $st=ts$ for certain pairs $s,t\in S$. We call $W$ a \emph{right-angled Coxeter group} if it admits at least one right-angled Coxeter generating set. The argument below shows, in particular, that if one Coxeter generating set for $W$ is right-angled, then every Coxeter generating set for $W$ is right-angled.

The main result is the following rigidity theorem.

\begin{theorem}\label{thm:main}
	Let $W$ be a right-angled Coxeter group. If $S$ and $S'$ are finite Coxeter generating sets for $W$, then the Coxeter systems $(W,S)$ and $(W,S')$ are equivalent.
\end{theorem}

The hypothesis that $W$ is right-angled is essential. For instance,
\[
	\langle a,b \mid a^2=b^2=(ab)^6=1\rangle
\]
and
\[
	\langle r,s,t \mid r^2=s^2=t^2=(rs)^3=(rt)^2=(st)^2=1\rangle
\]
are inequivalent Coxeter presentations of the symmetry group of a regular hexagon.

\section{The nerve of a Coxeter system}

Let $(W,S)$ be a Coxeter system. For $A\subseteq S$, let $W_A$ denote the subgroup generated by $A$, with $W_\varnothing=\{1\}$. The pair $(W_A,A)$ is again a Coxeter system \cite[p.~113]{Humphreys}, and standard parabolic subgroups satisfy
\[
	W_A\cap W_B=W_{A\cap B}
\]
for all $A,B\subseteq S$ \cite[p.~114]{Humphreys}. A subgroup of the form $W_A$ is called a \emph{standard parabolic subgroup}, and a conjugate of a standard parabolic subgroup is called a \emph{parabolic subgroup}.

The \emph{nerve} of $(W,S)$ is
\[
	\Nerve(W,S)=\{\,\varnothing\ne A\subseteq S : W_A \text{ is finite}\,\}.
\]
We also write
\[
	\MaxNerve(W,S)=\{\,A\in\Nerve(W,S): A\subseteq B\in\Nerve(W,S) \text{ implies } A=B\,\}
\]
for the collection of maximal simplices of the nerve.

If $(W,S)$ is right-angled, then the nerve is a flag complex: a nonempty subset $A\subseteq S$ lies in $\Nerve(W,S)$ if and only if every two-element subset of $A$ lies in $\Nerve(W,S)$. In this case, $W_A$ is finite precisely when the elements of $A$ commute pairwise, and then
\[
	W_A\cong (\Ztwo)^A.
\]

\section{Finite subgroups and the Davis--Vinberg complex}

We recall the standard description of finite subgroups of a Coxeter group. Let $(W,S)$ be a Coxeter system, and let $\Cosets$ be the set of all left cosets $wW_T$ for which $T\subseteq S$ and $W_T$ is finite. The \emph{Davis--Vinberg complex} $\Sigma(W,S)$ is the simplicial complex whose simplices are the nonempty finite chains in $\Cosets$, ordered by inclusion:
\[
	\Sigma(W,S)=\left\{\,F\subseteq\Cosets :
	\begin{array}{l}
		F\ne\varnothing,                                       \\
		\text{and any two elements of }F\text{ are comparable} \\
		\text{by inclusion}
	\end{array}
	\,\right\}.
\]
The group $W$ acts simplicially on $\Sigma(W,S)$ by left multiplication.

Let $X=|\Sigma(W,S)|$ be the geometric realization. There is a complete $W$-invariant metric on $X$ for which $X$ is a Hadamard space, that is, a complete simply connected geodesic metric space of nonpositive curvature in the sense of Alexandrov; see \cite{Moussong,Davis}. We use only the following consequence of this theory: every finite group of isometries of a Hadamard space has a nonempty fixed-point set \cite[Cor.~II.2.8]{BridsonHaefliger}.

\begin{lemma}\label{lem:finite-subgroups}
	Every finite subgroup of $W$ is contained in a finite parabolic subgroup. Consequently, the maximal finite subgroups of $W$ are exactly the groups
	\[
		wW_Aw^{-1},\qquad w\in W,
	\]
	where $A\in\MaxNerve(W,S)$.
\end{lemma}

\begin{proof}
	Let $G$ be a finite subgroup of $W$. By the fixed-point theorem quoted above, $G$ fixes a point of $X$. The action of $W$ on $\Sigma(W,S)$ is without inversions: if an element stabilizes a simplex, then it fixes each vertex of that simplex. Hence $G$ fixes a vertex of $\Sigma(W,S)$.

	The vertices of $\Sigma(W,S)$ are precisely the cosets $wW_A$ with $A\in\Nerve(W,S)\cup\{\varnothing\}$. If $G$ fixes the vertex $wW_A$, then $G$ is contained in its stabilizer, namely $wW_Aw^{-1}$. This is a finite parabolic subgroup. Maximal finite subgroups are therefore the maximal finite parabolic subgroups, which are exactly the conjugates $wW_Aw^{-1}$ with $A\in\MaxNerve(W,S)$.
\end{proof}

An alternative proof of Lemma~\ref{lem:finite-subgroups} is indicated in \cite[p.~130]{Bourbaki}.

\section{The right-angled case}

Let $(W,S)$ be right-angled, and let $S'$ be another finite Coxeter generating set for the same group $W$.

By Lemma~\ref{lem:finite-subgroups}, every finite subgroup of $W$ is contained in a finite parabolic subgroup for the system $(W,S)$. Since the system is right-angled, each finite parabolic subgroup is an elementary abelian $2$-group. Thus every element of $W$ has order $1$, $2$, or $\infty$. In any Coxeter system, $m(s',t')$ is the order of $s't'$. Hence $m(s',t')\in\{2,\infty\}$ for distinct $s',t'\in S'$, so $(W,S')$ is also right-angled.

Let $[W,W]$ denote the commutator subgroup, and let
\[
	q\colon W\longrightarrow W_{\ab}=W/[W,W]
\]
be the abelianization map. For a right-angled Coxeter system, the abelianization is isomorphic to $(\Ztwo)^S$. In particular, $|S|=|S'|$.

Let
\[
	V_S=\bigcup_{A\in\Nerve(W,S)\cup\{\varnothing\}} W_A
\]
be the union of the finite standard parabolic subgroups of $(W,S)$. Since $(W,S)$ is right-angled,
\[
	V_S=\{1\}\cup\{\,s_1s_2\cdots s_k : \{s_1,\ldots,s_k\}\in\Nerve(W,S),\ s_i\ne s_j\text{ for }i\ne j\,\}.
\]
The restriction of $q$ to $V_S$ is injective. Indeed, after identifying $W_{\ab}$ with $(\Ztwo)^S$,
we see that $q(s_1\cdots s_k)(t) = 1$ if $t \in \{s_1,\ldots,s_k\}\subseteq S$, and $0$ otherwise.

For each $A\in\MaxNerve(W,S)$, we shall associate a unique maximal simplex $A^*\in\MaxNerve(W,S')$.

\begin{theorem}\label{thm:maximal-correspondence}
	For every $A\in\MaxNerve(W,S)$, there exists a unique $A^*\in\MaxNerve(W,S')$ such that
	\[
		q(W_A)=q(W_{A^*}).
	\]
	Moreover, the map $A\mapsto A^*$ is a bijection from $\MaxNerve(W,S)$ to $\MaxNerve(W,S')$.
\end{theorem}

\begin{proof}
	Let $A\in\MaxNerve(W,S)$. Then $W_A$ is a maximal finite subgroup of $W$. Applying Lemma~\ref{lem:finite-subgroups} to the Coxeter system $(W,S')$, there are $w\in W$ and $A^*\in\MaxNerve(W,S')$ such that
	\[
		W_A=wW_{A^*}w^{-1}.
	\]
	Since $W_{\ab}$ is abelian, conjugate subgroups have the same image under $q$. Thus $q(W_A)=q(W_{A^*})$.

	The set $A^*$ is unique because $q$ is injective on the union $V_{S'}$ of the finite standard parabolic subgroups for $(W,S')$. The same argument with $S$ and $S'$ interchanged gives an inverse correspondence, so $A\mapsto A^*$ is a bijection.
\end{proof}

\begin{theorem}\label{thm:intersections}
	If $A_1,\ldots,A_r\in\MaxNerve(W,S)$, then
	\[
		\left|\bigcap_{i=1}^r A_i\right|=
		\left|\bigcap_{i=1}^r A_i^*\right|.
	\]
\end{theorem}

\begin{proof}
	Using the standard parabolic intersection formula and the injectivity of $q$ on $V_S$, we have
	\[
		\left|W_{\cap_i A_i}\right|
		=\left|\bigcap_{i=1}^r W_{A_i}\right|
		=\left|\bigcap_{i=1}^r q(W_{A_i})\right|.
	\]
	By Theorem~\ref{thm:maximal-correspondence}, $q(W_{A_i})=q(W_{A_i^*})$ for every $i$, and therefore
	\[
		\left|\bigcap_{i=1}^r q(W_{A_i})\right|
		=\left|\bigcap_{i=1}^r q(W_{A_i^*})\right|
		=\left|\bigcap_{i=1}^r W_{A_i^*}\right|
		=\left|W_{\cap_i A_i^*}\right|.
	\]
	Since all relevant standard parabolic subgroups are elementary abelian $2$-groups, $|W_B|=2^{|B|}$ for every simplex $B$ of either nerve. The result follows.
\end{proof}



\begin{proof}[Proof of Theorem~\ref{thm:main}]
	For each $s\in S$, define the \emph{membership pattern} of $s$ to be
	\[
		\mathcal A(s)=\{\,A\in\MaxNerve(W,S): s\in A\,\}.
	\]
	By inclusion--exclusion, the number of elements \(s\in S\) whose
	membership pattern is exactly \(\mathcal U\) is
	\[
		\sum_{\mathcal V\subseteq \MaxNerve(W,S)\setminus \mathcal U}
		(-1)^{|\mathcal V|}
		\left|
		\bigcap_{A\in \mathcal U\cup \mathcal V} A
		\right|.
	\]
	The corresponding number for \(S'\), with membership pattern
	\(\mathcal U^*=\{A^*:A\in\mathcal U\}\), is
	\[
		\sum_{\mathcal V\subseteq \MaxNerve(W,S)\setminus \mathcal U}
		(-1)^{|\mathcal V|}
		\left|
		\bigcap_{A\in \mathcal U\cup \mathcal V} A^*
		\right|.
	\]
	These two sums are equal term by term by Theorem~\ref{thm:intersections}.
	Hence the two systems have the same number of generators of each membership
	pattern.

	It follows that there is a bijection $\phi\colon S\to S'$ such that, for every $A\in\MaxNerve(W,S)$ and every $s\in S$,
	\[
		s\in A \quad\Longleftrightarrow\quad \phi(s)\in A^*.
	\]

	Now let $s,t\in S$ with $s\ne t$. Since $(W,S)$ is right-angled, $m(s,t)=2$ exactly when $\{s,t\}$ is contained in some maximal simplex $A\in\MaxNerve(W,S)$. By the defining property of $\phi$, this occurs exactly when $\{\phi(s),\phi(t)\}$ is contained in the corresponding maximal simplex $A^*\in\MaxNerve(W,S')$, which is exactly the condition that $m(\phi(s),\phi(t))=2$. Otherwise both exponents are $\infty$.

	Thus
	\[
		m(s,t)=m(\phi(s),\phi(t))
	\]
	for all $s,t\in S$. Therefore $\phi$ respects the Coxeter matrices and extends to an automorphism of $W$ carrying $S$ onto $S'$. Hence $(W,S)$ and $(W,S')$ are equivalent.
\end{proof}

\end{document}